
\input amstex 
\documentstyle{amsppt}
\input bull-ppt
\keyedby{bull317e/kmt}

\topmatter
\cvol{27}
\cvolyear{1992}
\cmonth{October}
\cyear{1992}
\cvolno{2}
\cpgs{273-278}
\title Additive functions on shifted primes \endtitle
\author P. D. T. A. Elliott\endauthor
\address Department of Mathematics, University of Colorado, 
Boulder, Colorado 80309-0001\endaddress
\date November 26, 1991. Presented at the 1992 Illinois 
Number Theory
Conference, University of Illinois at Urbana-Champaign, 
April 3--4, 1992\enddate
\abstract Best possible bounds are obtained for the
concentration function of an additive arithmetic function 
on sequences
of shifted primes.\endabstract
\endtopmatter

\document



A real-valued function $f$ defined on the positive 
integers is {\it
additive} ~if it satisfies $f(rs)  = f(r)+f(s)$ whenever 
$r$ and $s$ are
coprime.  Such functions are determined by their values on 
the
prime-powers.

For additive arithmetic function $f$, let $C_h$ denote
the frequency amongst the integers $n$ not exceeding $x$ 
of those for
which $h < f(n) \leq h+1$.  Estimates for $C_h$ that are 
uniform in $h,
f$, and $x$ play a vital r\^ole in the study of the value 
distribution of
additive functions.  They can be employed to develop 
criteria necessary
and sufficient that a suitably renormalised additive 
function possess a
limiting distribution, as well as to elucidate the 
resulting limit law.
They bear upon problems of algebraic nature, such as the 
product and
quotient representation of rationals by rationals of a 
given type.  In
that context their quantitative aspect is important.

It is convenient to write $a \ll b$ {\it uniformly in} 
~$\alpha$ if on
the values of $\alpha$ being considered the functions  $a, 
b$ satisfy
$|a(\alpha )| \leq cb(\alpha )$ for some absolute constant 
$c$.  When
the uniformity is clear, I do not declare it.

Let
$$W(x) = 4 + \min_\lambda \left( \lambda^2 + \sum_{p \leq 
x} \frac{1}{p}
~\min (1, |f(p) - \lambda \log p|)^2 \right) ,$$
where the sum is taken over the prime numbers.
Improving upon an earlier result of Hal\'asz, Ruzsa proved 
that $C_h \ll
W(x)^{-1/2}$, uniformly in $h, f$, and  $x \geq 2$ [9].  
This result is 
best possible in the sense that for
each of a wide class of additive functions there is a 
value of $h$ so
that the inequality goes the other way.

From a number theoretical point of view it is desirable to 
possess
analogs of Ruzsa's result in which the additive function 
$f$ is confined
to a particular sequence of integers of arithmetic 
interest.  In this
announcement I consider shifted primes.

Let $a$ be a nonzero integer.  Let $Q_h$ denote the 
frequency amongst
the primes $p$ not exceeding $x$ of those for which $h < 
f(p+a) \leq
h+1$.

\proclaim{Theorem 1}  The estimate  $Q_h \ll W(x)^{-1/2}$ 
holds
uniformly in $h, f$, and $x \geq 2$.
\endproclaim

If for an integer $N \geq 3$ we define $S_h$ to be the 
frequency amongst
the primes $p$ less than $N$ of those for which $h < 
f(N-p) \leq h+1$,
and set 
$$Y(N) = 4 + \min_\lambda \left( \lambda^2 + \sum \Sb p < 
N \\
(p,N)=1\endSb \frac{1}{p} ~\min (1, |f(p)-\lambda \log 
p|)^2 \right) ,$$
then there is an analogous result.

\proclaim{Theorem 2}  The estimate $S_h \ll Y(N)^{-1/2}$ 
holds uniformly
in $h, f$, and $N \geq 3$.
\endproclaim

The estimates given in these two theorems are of the same 
quality as
Ruzsa's and again best possible.  In particular, Theorem 1 
improves the
bound $Q_h \ll W(x)^{-1/2} (\log W (x))^2$ of Timofeev 
[10].  If
$$E(x) = 4 + \sum \Sb p \leq x \\ f(p) \not= 0\endSb 
\frac{1}{p} ,$$
then Timofeev shows that the number of primes not 
exceeding $x$ for
which $f(p+a)$ assumes any (particular) value is $\ll \pi 
(x)E(x)^{-1/2}
(\log E(x))^2$.  Employing the present Theorem 1, the 
logarithmic factor
may be stripped from this bound.  The improved inequality 
is then
analogous to an estimate of Hal\'asz concerning additive 
functions on
the natural numbers and, in a sense, best possible [8].

The concentration function estimate of Theorem 2 also has 
many
applications, in particular, to the study of the value 
distribution of
additive functions.  These are new and of a new type.  
They involve not
only the primes but also the length of the interval on 
which the
additive function is considered.  Thus the frequencies
$$(\pi (N-1))^{-1} \sum \Sb p < N \\ f(N-p) \leq z\endSb 1$$
possess a limiting distribution function as $N \rightarrow 
\infty$ if
and only if the three series
$$\sum_{|f(p)| > 1} \frac{1}{p} , \qquad \sum_{|f(p)| \leq 
 1}
\frac{f(p)}{p} , \qquad \sum_{|f(p)| \leq 1} 
\frac{f(p)^2}{p}$$
converge.  The latter is the classical condition of 
Erd\"os and Wintner
required when considering frequencies over the natural 
numbers [7].
More complicated examples involving unbounded 
renormalisations of
additive functions can also be successfully treated.

The method of this paper lends itself well to the study of 
the
representation of rationals by products and quotients of 
shifted primes.

The proofs of Theorems 1 and 2 apply Fourier analysis.  
Since the
F\'ejer kernel is nonnegative, $Q_h$ does not exceed
$$3\pi (x)^{-1} \sum_{p \leq x} \int_{-1}^1 
(1-|t|)e^{-ith} g(p+a)\,dt,
\tag 1$$
where $g(n)$ is the multiplicative function $\exp 
(itf(n))$.  To deal
directly with the mean value of $g(p+a)$ over the primes 
would require
finer information concerning the distribution of primes in 
residue
classes than is currently available.  Let $3|a| \leq w 
\leq z$.
Ultimately $z$ will be chosen a power of $x$, $w$ a power 
of $\log x$.
Let $P, R$ denote the products of the primes in the ranges 
$3|a| < p
\leq w$, $w < p \leq z$ respectively.  I majorize (1) by 
introducing
a Selberg square function $(\sum_{d\mid (n,R)}
\lambda_d)^2$, where the $\lambda_d$ are real, zero if $d 
> z$,
$\lambda_1 = 1$.  Expanding and interchanging the order of 
summation
gives
$$Q_h \leq \frac{3}{\pi (x)} \int_{-1}^1 (1-|t|)e^{-ith} 
\sum_{d_j \mid
R} \lambda_{d_1} \lambda_{d_2} \sum \Sb n \leq x, 
(n,P)=1\\ n \equiv 0
~(\text{mod} \,[d_1 ,d_2 ])\endSb g(n+a)\,dt + 
\frac{3z}{\pi (x)} .\tag
2$$
We are reduced to the study of multiplicative functions on 
arithmetic
progressions with moduli large compared to $x$.  It may 
seem curious to
retain the condition $(n,P)=1$.  However, the choice of a 
nonprincipal
character (mod 3) for $g$ shows that the expected estimate
$$\sum \Sb n \leq x \\ n \equiv r ~(\text{mod} \,D)\endSb 
g(n) =
\frac{1}{\phi (D)} \sum \Sb n \leq x \\ (n,D)=1\endSb g(n) +
\text{`small'} $$
is in general false.  In [1, 4, 6] it is shown that the 
moduli
$D$ for which such an estimate fails to be reasonably true 
are multiples
of a single modulus $D_0$.  The present situation is 
arranged so that
the complications due to the existence of $D_0$ are bound 
up in the
condition $(n,P)=1$ and that effectively $D_0 ,R$ have no 
common
divisors.

The moduli $d_j$ dividing $R$, with $d_j \leq z$ are dealt 
with by means
of the following result.

For a multiplicative function $g$, with values in the 
complex unit
disc, define an exponentially multiplicative
function $g_1$ by $g_1 (p^k ) = g(p)^k /k!$.  Define the 
multiplicative
function $h$ by convolution: ~$g = h\ast g_1$.  Thus $g(p) 
= g_1 (p)$,
$h(p) = 0$.  Moreover, $|h(p^k )| \leq e$.  For $A \geq 0$ 
define
$$\align
\beta_1 (n)& = \sum \Sb ump = n \\ u \leq (\log x)^{2A} \\ 
p \leq (\log
x)^{6A+15} \endSb h(u) g_1 (m) g(p) \frac{\log p}{\log mp} 
,\\
\beta_2 (n)& = \sum \Sb urp = n \\ u \leq (\log x)^{2A} \\ 
r \leq (\log
x)^{6A+15} \endSb h(u) g_1 (r) g(p) \frac{\log p}{\log rp} 
,\endalign
$$
and set $\beta (n) = g(n) - \beta_1 (n) - \beta_2 (n)$.  
Note that
$\beta_j (n) \ll 1$ uniformly in $n, j$.

\proclaim{Lemma 1}  Let $0 < \delta < 1/2$.  Then
$$\align
&\sum_{D_1 D_2 \leq x^\delta} \max_{(r,D_1 D_2 )=1} 
\max_{y \leq x} \left|
\sum \Sb n \leq y \\ n \equiv r ~(\text{mod} \,D_1 D_2 
)\endSb \beta (n)
- \frac{1}{\phi (D)} \sum \Sb n \leq y, (n,D_2 )=1 \\ n 
\equiv r
~(\text{mod} \,D_1 )\endSb \beta (n) \right|\\
&\qquad\ll x(\log x)^{-A} (\log \log x)^2 + w^{-1} x(\log 
x)^{2A+8} (\log\log
x)^2\\
&\qquad\quad+ w^{-1/2} x(\log x)^{5/2} \log\log x,\endalign
$$
where $D_1$ is confined to integers whose prime factors do 
not exceed
$w$ and $D_2$ to integers all of whose prime factors 
exceed $w$.  The
implied constant depends at most upon $\delta ,A$.
\endproclaim

Lemma 1 represents a generalisation to largely arbitrary 
multiplicative
functions of the well-known theorem of Bombieri and 
Vinogradov
concerning primes in arithmetic progressions.  The 
parameter $\delta$
may be replaced by $1/2 - \varepsilon (x)$ for a certain 
positive
function $\varepsilon (x)$, which approaches zero as $x 
\rightarrow
\infty$.  Of importance here is the quality of the error 
term.  For $w
\geq (\log x)^{3A+8}$ it is as good as that of Bombieri 
and Vinogradov.
To this end the functions $\beta_j$ were introduced, 
manifesting the
assertion of [5, p. 408], already in view in [3, p. 178], 
that for
general multiplicative functions a change of form would be 
required.  In
particular, $\beta_2 (n)$ is largely supported on the 
primes and cannot
be removed without further information concerning $g$.  
Most integers
$n$ will have few prime divisors, so that effectively the 
$\beta_j (n)$
are $\ll \log\log x/\log x$ over the range $2 \leq n \leq 
x$.

The functions $\beta_j$ run through the treatment of the 
integral at (2)
along with the central function $g$.  A notable feature of 
the method is
the casting of the Selberg square functions on the 
multiplicative
integers in a r\^ole, which on the additive group of 
reals, is
traditionally played by a F\'ejer kernel.  The outcome is 
the estimate
$$Q_h \ll x^{-1} \log w \int_{-1}^1 (1-|t|)e^{-ith} \sum 
\Sb n \leq x
\\ (n, P) = 1\endSb g(n+a)\,dt + (\log x)^{-1} (\log\log 
x)^2 .\tag 4$$
The complications introduced by the exceptional modulus 
$D_0$ mentioned
earlier must now be dealt with.  To this end [4] or [6] 
may be applied.
For simplicity of exposition I appeal to Theorem 1 of [6].

\proclaim{Lemma 2}  Let $0 < \gamma < 1$, $0 < \delta < 
1/8$, $2 \leq
\log N \leq Q \leq N$.  Then any multiplicative function 
$g$ with
values in the complex unit disc satisfies
$$\sum \Sb n \leq x \\ n \equiv r ~(\text{mod} \,D)\endSb 
g(n) =
\frac{1}{\phi (D)} \sum \Sb n \leq x \\ (n,D)=1\endSb g(n) 
+ O\left(
\frac{x}{\phi (D)} \Big( \frac{\log Q}{\log x} \Big)^{1/8 
-\delta}
\right)$$
uniformly for $N^\gamma \leq x \leq N$, for all $(r,D)=1$, 
for all $D
\leq Q$ save possibly for the multiples of a $D_0 > 1$.
\endproclaim

From Lemma 2 with $N=x$, $Q = \exp ((\log\log x)^2 )$ I 
obtain the
following estimate.

\proclaim{Lemma 3}  Let $w$ be a power of $\log x$ and $P$ 
the product of
the primes in the interval $(y, w]$, where $3|a| \leq y 
\leq w$.  Then
either
$$\sum \Sb n \leq x \\ (n-a, P)=1\endSb g(n) = \prod_{y < 
p \leq w}
\left( 1 - \frac{1}{p-1} \right) \sum_{n \leq x} g(n) 
\prod \Sb p \mid n
\\ p > y\endSb \left( \frac{p-1}{p-2} \right) + O(x(\log 
x)^{-1/10} ),$$
or there is a prime divisor $q$ of $P$ such that
$$\align\sum \Sb n \leq x \\(n-a, P)=1\endSb g(n) = &\prod 
\Sb  y < p \leq w
\\ p \not= q\endSb \left( 1- \frac{1}{p-1} \right) \sum 
\Sb n \leq x \\
(n,q)=1\endSb g(n) \prod \Sb p \mid n \\ p > y \endSb \left(
\frac{p-1}{p-2} \right)\\
& + O(x(\log x)^{-1/10} ).\endalign$$
The implied constants do not depend upon $y, g$, or $q$.
\endproclaim

The prime $q$ may vary with $g$ and $x$.

It follows from (4) and Lemma 3 that
$$Q_h \ll x^{-1} \int _{-1}^1 (1-|t|)e^{-ith} \sum_{n \leq 
x} g(n) \prod
\Sb p \mid n \\ p > 3|a| \endSb \left( \frac{p-1}{p-2} 
\right) dt +
(\log x)^{-1/12} ,\tag 5$$
with possibly a condition $(n,q)=1$ required in the sum.  
Whilst the
function $g$ in Lemmas 1, 2, and 3 may be arbitrary up to 
having values in the
unit complex disc, in (5) $g$ has the special form $\exp 
(itf(n))$.  The
exceptional prime $q$ may therefore vary with  $t$.  It 
can be arranged
that $q$ may only exist on intervals, on each of which it 
will be
constant.  The integral at (5) is therefore well defined.  
Without the
condition $(n,q)=1$ we may now follow the original 
treatment of Ruzsa
[9], who considered a similar integral without the weight 
factor $\Pi
(p-1)/(p-2)$.  The extra condition $(n,q)=1$ introduces 
some further
complications, but they can be overcome.  

Similarly
$$S_h \ll \phi (N)^{-1} \int_{-1}^1 (1-|t|)e^{-ith} \sum 
\Sb n \leq N
\\ (n,N)=1\endSb g(n) \prod \Sb p \mid n \\ p > 3\endSb 
\left(
\frac{p-1}{p-2} \right) dt + (\log N)^{-1/10} .$$
Once again an auxiliary condition $(n,q)=1$ may be needed 
in the sum.
Since $N$ may have many prime factors, the condition 
$(n,N)=1$
introduces a new complication, but this, too, can be 
overcome.  It may
be remarked here that Theorem 2 of [2] shows that in quite 
general
circumstances conditions of the type $(n,N)=1$ may be 
factored out of
mean values of multiplicative functions.

It transpires that the parameter $\lambda$ appearing in 
the definitions
of $W(x)$ and $Y(N)$ may be restricted by $|\lambda | \leq 
(\log x)^2$,
$|\lambda | \leq (\log N)^2$ respectively.


\enddocument